\documentclass[11pt]{article}
\usepackage[margin=1in]{geometry} 
\usepackage{enumerate}
\usepackage{amsmath,amsthm,amssymb,amsfonts, fancyhdr, comment, parskip, tikz, caption, subcaption, graphicx,tikz-network}
\usepackage{diagbox}
\usepackage[colorlinks]{hyperref}
\usetikzlibrary{graphs}
\usetikzlibrary{shadings}
\usetikzlibrary{shapes.geometric}
\usetikzlibrary{arrows}
\usetikzlibrary{calc, positioning, fit, shapes.misc}

\usepackage{todonotes}

\definecolor{purple}{RGB}{148,0,211}
\definecolor{green}{RGB}{50,205,50}

\newtheorem{theorem}{Theorem}[section]

\newtheorem{lemma}[theorem]{Lemma}

\newtheorem{proposition}[theorem]{Proposition}

\theoremstyle{definition}
\newtheorem{example}[theorem]{Example}
\newtheorem{definition}[theorem]{Definition}

\newcommand{\mc}[1]{\mathcal{#1}}
\newcommand{\cbrac}[1]{\left\{ #1\right\}} 
\newcommand{\floor}[1]{\lfloor#1\rfloor}
\tikzstyle{arch} = [out=30, in=150]

\newcommand{\R}{\textup{R}}
\newcommand{\AR}{\reflectbox{\textup{R}}}

\providecommand{\MR}[1]{}

\DeclareMathOperator{\ex}{ex}

\title{Generalized anti-Ramsey numbers}
\author{Deepak Bal \footnote{Department of Mathematics, Montclair State University, Montclair, NJ. \texttt{deepak.bal@montclair.edu}.} \and Jonathan Cutler\footnote{Department of Mathematics, Montclair State University, Montclair, NJ. \texttt{jonathan.cutler@montclair.edu}.}}

\date{}

\begin{document}

\maketitle
\begin{abstract}
    The anti-Ramsey number was introduced by Erd\H{o}s, Simonovits, and S\'os and is defined to be the maximum number of colors that can be used to color the edges of a graph to avoid a rainbow copy of some graph.  One can generalize this by, instead of forbidding a rainbow copy of a graph, putting a cap on the number of colors that appear on any copy of that graph.  To be precise, let $\AR(G,H,q)$ be the maximum number of colors that can be used to color $E(G)$ so that every copy of $H$ sees at most $q$ colors.  One can see that $\AR(G,H,q)$ is nondecreasing in $q$.  In this paper, we investigate the linear and quadratic thresholds of $\AR(G,H,q)$ when $G=K_n$ and $H$ is connected.
\end{abstract}
\section{Introduction}

Ramsey theory is of course a topic that appears in many areas of mathematics.  The original questions investigated by Ramsey involved finding the minimum size of a structure that guaranteed the existence of a monochromatic substructure.  Erd\H{o}s, Simonovits, and S\'os \cite{ESS} introduced anti-Ramsey problems which ask for the maximum number of colors that can be used to color a structure to avoid a rainbow copy of some substructure.  These concepts have been generalized in numerous ways since their introduction and this paper will investigate one such generalization.

\begin{definition}
The \emph{generalized Ramsey number} $\R(G,H,q)$ is the minimum number of colors in a coloring of $E(G)$ such that every copy of $H$ receives at least $q$ colors.   We may analogously define the \emph{generalized anti-Ramsey number} $\AR(G,H,q)$ as the maximum number of colors in a coloring of $E(G)$ such that every copy of $H$ receives at most $q$ colors.\footnote{We would like to thank Peter Winkler for suggesting this notation.} 
\end{definition}

Note that for $\R(G,H,q)$, the case $q=2$ corresponds to the classical multicolor Ramsey problem. For $\AR(G,H,q)$, the case $q = e(H)-1$ corresponds to the classical anti-Ramsey problem.  It should be noted that Mubayi and West~\cite{MW} studied an even more general extension of the anti-Ramsey number where the number of colors in any copy of $H$ is bounded both from above and below.

These definitions can take a bit to get used to and so, in an effort to make them clearer, we note the following.  Call an $M$-coloring, $c:E(G)\to [M]$, \emph{exact} if it is surjective, i.e., if each color appears at least once. In this case, we say that $E(G)$ has been \emph{exactly colored} with $M$ colors.  Then:

\begin{itemize}
\item $\R(G,H,q)\le M$ means that there exists an $M$-coloring of $E(G)$ such that every copy of $H$ receives at least $q$ colors.

\item $\R(G,H,q)> M$ means that every $M$-coloring of $E(G)$ contains a copy of $H$ with at most $q-1$ colors.

\item $\AR(G,H,q) < M$ means that every exact $M$-coloring of $E(G)$  contains a copy of $H$ with at least $q+1$ colors.

\item $\AR(G,H,q) \ge M$ means that there exists an exact $M$-coloring of $E(G)$ such that every copy of $H$ receives at most $q$ colors.
\end{itemize}


We give a short and incomplete history of generalized Ramsey numbers before moving on to anti-Ramsey problems.  The problem of determining $\R(K_n,K_k,q)$ was first introduced by Erd\H{o}s and Shelah (see \cite{E64}) and the first results were found by Elekes, Erd\H{o}s, and F\"uredi (see \cite{E74}).  Erd\H{o}s and Gy\'arf\'as \cite{EG97} then found the linear and quadratic thresholds for $\R(K_n,K_k,q)$, i.e., the smallest $q$ for which $\R(K_n,K_k,q)=\Omega(n)$ and $\Omega(n^2)$, respectively.  Axenovich, F\"uredi, and Mubayi~\cite{AFM} investigated thresholds for $\R(K_{n,n},K_{k,k},q)$ and, more recently, Krueger~\cite{K20} studied thresholds for $\R(K_n,P_k,q)$. The first author, along with Bennett, Heath, and Zerbib~\cite{BBHZ} generalized the problem further to requiring that each graph in some set sees at least $q$ colors, and also to the hypergraph setting. Throughout the paper, unless otherwise stated, assume that target graphs $H$ are fixed in asymptotic statements and that $n$ is sufficiently large.

The story of generalized anti-Ramsey numbers begins with Erd\H{o}s, Simonovits, and S\'os \cite{ESS} who introduced the classical anti-Ramsey number, $\AR(G,H,e(H)-1)$.  They conjectured that 
\[
    \AR(K_n,C_k,k-1)=n\left(\frac{k-2}2+\frac{1}{k-1}\right)+O(1),
\]
and proved the $k=3$ case.  Alon \cite{A} proved the $k=4$ case, and Jiang and West~\cite{JW} improved the general lower bound, before Montellano-Ballesteros and Neumann-Lara \cite{MBNL} resolved the full conjecture.

\begin{theorem}[Montellano-Ballesteros, Neumann-Lara]\label{thm:mbnl}
    For all $n\geq k\geq 3$, we have
    \[
        \AR(K_n,C_k,k-1)=n\left(\frac{k-2}2+\frac{1}{k-1}\right)+O_k(1).
    \]
\end{theorem}

Jiang~\cite{J} studied the problem when $H$ is a star.  We let $S_{r+1}$ be the star with $r+1$ leaves, i.e., $K_{1,r+1}$.  Jiang proved that if $r\leq n-2$, then
\[
    \AR(K_n,K_{1,r+1},r)=\left\lfloor\frac{1}2n(r-1)\right\rfloor+\left\lfloor \frac{n}{n-r+1}\right\rfloor +\epsilon,
\]
where $\epsilon\in\{0,1\}$.  Note that this also gives the generalized anti-Ramsey number for the star: every $K_{1,r}$ sees at most $q$ colors is equivalent to containing no rainbow $K_{1,q+1}$.   Jiang and West~\cite{JW} investigated the anti-Ramsey problem where $H$ is a tree or indeed a family of trees.  They showed that if $T$ is a tree with $k$ edges and $n\geq 2k$, then
\[
    \frac{n}2\left\lfloor \frac{k-2}2\right\rfloor +c_k\leq \AR(K_n,T,k-1)\leq n(k-1),
\]
where $c_k$ is a constant depending only on $k$.

As with generalized Ramsey numbers, we will be interested in the behavior of $\AR(K_n,G,q)$ for a fixed graph $G$ as $q$ increases.  Note that
\begin{enumerate}
    \item $\AR(K_n,G,1)=1$, provided that $e(G)\geq 2$ and $n\geq n(G)$, and
    \item $\AR(K_n,G,e(G))=\binom{n}2$,
\end{enumerate}
since in the first case, one can use at most one color on the edges of $K_n$, and in the second, every edge of $K_n$ can receive a different color.  So, for a given graph $G$, we will investigate the smallest integer $q$ for which $\AR(K_n,G,q) = \Omega(n)$, which we will call the \emph{linear threshold}.  Similarly, we will be interested in finding the \emph{quadratic threshold} as well. 

Perhaps the most studied question in this vein is that of determining $\AR(K_n,K_p,q)$ for various values of $n$, $p$, and $q$.  We begin with a construction that gives a lower bound.

\begin{example}\label{ex:star}
    Suppose the vertices of $K_n$ are labeled $\{1,2,\ldots,n\}$.  We get an edge-coloring of $K_n$ by letting edge $ij$ for $i<j$ be colored $i$.  Then note that any $K_p$ for $1\leq p\leq n$ sees exactly $p-1$ colors since the edges are colored with colors corresponding to the vertex labels except for the largest one.
\end{example}

Note that Example~\ref{ex:star} gives that $\AR(K_n,K_p,p-1)\geq n-1$.  Erd\H{o}s, Simonovits, and S\'os~\cite{ESS} proved that equality holds for $p=3$.  Ahlswede, Cai, and Zhang~\cite{ACZ} then showed equality holds for all $p$.
In fact, Ahlswede, Cai, and Zhang gave a result that covers other values of $q$ as well.  We let $t(n,k)$ be the Tur\'an number, so that it is the number of edges in the Tur\'an graph, $T(n,k)$.\footnote{The Tur\'an graph $T(n,k)$ is the complete $k$-partite graph on $n$ vertices whose part sizes differ by at most one, so that each part has size $\lceil n/k\rceil$ or $\lfloor n/k\rfloor$.}

\begin{theorem}[Ahlswede, Cai, Zhang]
    If $n\geq p\geq 3, k\ge 1$ and $1\le q' \le \floor{\floor{p/k}/2}$ , we have
    \[
        \AR(K_n,K_p,t(p,k)+q')=t(n,k)+q'.
    \]
\end{theorem}

Very roughly speaking, this result states that if $k\ge 2$, the number of colors that a $K_p$ sees is roughly quadratic in $p$ (about $(1-1/k)p^2/2$), then the generalized anti-Ramsey number is quadratic in $n$.  

Axenovich and K\"undgen~\cite{AK} investigated a generalization of this problem and noted that $\AR(K_n,K_p,q)=O(n)$ if and only if $1\leq q\leq p-1$.  In so doing, they noted that there are constructions of graphs with girth $p+1$ and at least $cn^{1+\varepsilon}$ edges.  (See, e.g., \cite{LUW}.)  If one colors the edges of such a graph with all different colors and its complement with one additional color, one gets a coloring with a super-linear number of colors in which every $K_p$ sees at most $p$ colors.  In other words, we get
\[
    \AR(K_n,K_p,p)\geq cn^{1+\varepsilon}.
\]

In the case $p=3$, all cases are covered by the results outlined above.  In particular, we have $\AR(K_n,K_3,2)=n-1$ by the result of Erd\H{o}s, Simonovits, and S\'os as of course $C_3=K_3$.  In the case $p=4$, we know $\AR(K_n,K_4,3)=n-1$.  However, the case $q=4$ is not covered by any of these results.  
It has been shown that 
\[
    \AR(K_n,K_4,4)=\Theta(n^{3/2}),
\]
where the lower bound comes from a result of Erd\H{o}s~\cite{E} that there is a graph $G$ on $n$ vertices with girth five and at least $cn^{3/2}$ edges.  The upper bound follows from a result of Montellano-Ballesteros~\cite{MB} on the anti-Ramsey number of diamonds.

In this paper, we investigate several problems related to the generalized anti-Ramsey number $\AR(G,H,q)$.  In Section~\ref{sec:pathcyc}, we present some results when $G$ is complete and $H$ is a path or a cycle.   In Section~\ref{sec:linear}, we prove our main result, in which we determine the linear threshold for connected graphs.  We write $\Delta(G)$ for the maximum degree of a graph $G$ and $\nu(G)$ for the maximum size of a matching in $G$.

\begin{theorem}\label{thm:trees}
    Let $G$ be a connected graph on at least 3 vertices and suppose $k=\min\cbrac{\nu(G), \Delta(G)}$. Assuming $n$ is sufficiently large, we have 
    \begin{enumerate}[(i)]
        \item\label{thm:tree_const} $\AR(K_n, G, k)\le (k+2)(k+1)$, and
        \item\label{thm:tree_linear}$\AR(K_n, G, k+1) =\Omega(n)$.
    \end{enumerate}
\end{theorem}

In Section~\ref{sec:quadratic} we prove the following theorem which determines the quadratic threshold for any graph $G$, which depends on whether $G$ is bipartite or not.

\begin{theorem}\label{thm:quad}
    Let $G$ be a graph on $v$ vertices with $e(G)\ge 2$ and suppose that $k$ is the number of edges in the largest bipartite subgraph of $G$. Then
    \begin{enumerate}[(i)]
        \item if $G$ is bipartite, then $k=e(G)$ and  $\AR(K_n, G, k-1) = O(n^{2-2/v})$ and $\AR(K_n, G, k) = \binom{n}2$, and
        \item if $G$ is nonbipartite, then $\AR(K_n, G, k) = O(n^{2-1/(v+1)})$ and $\AR(K_n, G, k+1) = \Omega(n^2)$.
    \end{enumerate}

\end{theorem}


\section{Paths and Cycles}\label{sec:pathcyc}



In this section, we investigate and gather results about $\AR(K_n,P_k,q)$ and $\AR(K_n,C_k,q)$ for various $n$, $k$, and $q$.  Much of this will follow from classical results, but we believe there is some value in collecting it for reference.  Throughout this section, we'll be interested in determining the asymptotics of $\AR(K_n,P_k,q)$ and $\AR(K_n,C_k,q)$ as $n$ goes to infinity.

We'll start with paths.  It is trivial that $\AR(K_n,P_k,k-1)=\AR(K_n,P_k,e(P_k))=\binom{n}2$ as there is no restriction on the coloring in this case.  But note that $\AR(K_n,P_k,k-2)$ corresponds exactly to the classical anti-Ramsey question as by stipulating that each $P_k$ sees at most $k-2$ colors is equivalent to containing no rainbow $P_k$s.

Erd\H{o}s, Simonovits, and S\'os~\cite{ESS} proved the following.
\begin{theorem}[Erd\H{o}s, Simonovits, S\'os]\label{thm:ess}
    Let $t\ge 0$ be a given integer, $\varepsilon\in \{0,1\}$, and $k=2t+3+\varepsilon$.  Then, if $n$ is large, we have
    \[
        \AR(K_n,P_k,k-2)=tn-\binom{t+1}2+1+\varepsilon.
    \]
\end{theorem}
The proof appeared in a paper of Simonovits and S\'os~\cite{SS}.  Further, Erd\H{o}s, Simonovits, and S\'os made a conjecture about the behavior for all $n$.  The full conjecture has been resolved by Yuan~\cite{Y}.  

Note that, for any graph $G$, we have that $\AR(K_n,G,q)$ is monotone in subgraphs.  To be precise, we have the following.

\begin{proposition}\label{prop:subgr}
    If $H$ is a subgraph of $G$ and $n\ge |V(G)|$, then 
    \[
        \AR(K_n,G,q)\leq \AR(K_n,H,q).
    \]
\end{proposition}
\begin{proof} 
To see this, observe that if there exists an $M$-coloring of $K_n$ such that every copy of $G$ receives at most $q$ colors, then certainly every copy of $H$ receives at most $q$ colors.  
\end{proof}
But then for all $\ell\geq k$, by Proposition~\ref{prop:subgr}, since $P_{\ell}\supseteq P_k$, we know $\AR(K_n,P_{\ell},k-2)=\AR(K_n,P_k,k-2)$ as the construction that works for smaller paths also works for longer paths.  Thus, we can give the asymptotics of $\AR(K_n,P_k,q)$ for all values of $k$ and $q$ by combining Theorem~\ref{thm:ess} and Proposition~\ref{prop:subgr}.  Note that we assume that $q\leq k-2$ since if $q=k-1$, we trivially have $\AR(K_n,P_k,k-1)=\binom{n}2$.
\begin{theorem}\label{thm:paths}
    If $n$ is sufficiently large, $1\le q\leq k-2$, and $q=2t+1+\varepsilon$ where $\varepsilon\in \{0,1\}$, then 
    \[
        \AR(K_n,P_k,q)=\AR(K_n,P_{q+2},q)=tn-\binom{t+1}2+1+\varepsilon.
    \]
\end{theorem}

\begin{proof}
    As noted above, we have that $\AR(K_n,P_k,q)\le\AR(K_n,P_{q+2},q)$ by Proposition~\ref{prop:subgr}.  Then, in Theorem~\ref{thm:ess}, we have $k=2t+3+\varepsilon$ so that if $k=q+2$, we know $q=2t+1+\varepsilon$, and so the second equality holds.  
    We end by showing a construction that gives the result.  This is a construction due to Erd\H{o}s, Simonovits, and S\'os~\cite{ESS}.   Choose $t$ vertices in $K_n$, say $x_1,x_2,\ldots,x_t$.  Color all edges with at least one endpoint an $x_i$ for some $i$ with different colors and color all other edges (so those in $V(K_n)\setminus \{x_1,\ldots,x_t\}$) with either one additional color, if $\varepsilon=0$, or two additional colors, if $\varepsilon=1$.  Then any $P_k$ sees at most $2t+1+\varepsilon=q$ colors.
\end{proof}

Moving on to cycles, since $e(C_k)=k$, we know that $\AR(K_n,C_k,k)=\binom{n}2$ for any $k\geq 3$.  We also know that $\AR(K_n,C_k,k-1)$ corresponds to the classical anti-Ramsey problem and so these are done by the theorem of Montellano-Ballesteros and Neumann-Lara~\cite{MBNL}, that is, Theorem~\ref{thm:mbnl}.  For other values of $q$, we can fall back on what we know about paths, as paths are subgraphs of cycles.

\begin{proposition}
    If $\ell\geq k\ge 3$ and $n$ is sufficiently large, then
    \[
        \AR(K_n,C_{\ell},k-2)=\AR(K_n,C_k,k-2)=\AR(K_n,P_k,k-2).
    \]
\end{proposition}

\begin{proof}
    For $\ell\geq k$, we know that $\AR(K_n,C_{\ell},k-2)\leq \AR(K_n,P_k,k-2)$ and $\AR(K_n,C_{k},k-2)\leq \AR(K_n,P_k,k-2)$ by Proposition~\ref{prop:subgr}.  So, we need to show equality in both cases.  We note that the same construction that works for the anti-Ramsey case of the path on $k$ vertices works for any cycle of length at least $k$.  (See the proof of Theorem~\ref{thm:paths}.) 
    So, we have equality across the board and the proposition holds.
\end{proof}

We summarize the generalized anti-Ramsey numbers of some small cycles in Table~\ref{tab:cycles}.

\renewcommand{\arraystretch}{1.5}
\begin{table}
\begin{center}
\begin{tabular}{l||c|c|c|c|c|c|c|c}
\diagbox{$k$}{$q$} & 1 & 2 & 3 & 4 & 5 & 6 & 7 &8\\ \hline\hline
3 & 1 & $n-1$ & $\binom{n}2$ & & & & &\\ \hline
4 & 1 & 2 & $\frac{4}3n+O(1)$ & $\binom{n}2$ & & & &\\ \hline
5 & 1 & 2 & $n$ & $\frac{7}4n+O(1)$& $\binom{n}2$ & & & \\ \hline
6 & 1 & 2 & $n$ & $n+1$ & $\frac{11}5n+O(1)$ & $\binom{n}2$ & &\\ \hline
7 & 1 & 2 & $n$ & $n+1$ & $2n-2$ & $\frac{8}3n+O(1)$ & $\binom{n}2$ & \\ \hline
8 & 1 & 2 & $n$ & $n+1$ & $2n-2$ & $2n-1$ & $\frac{22}7n+O(1)$ & $\binom{n}2$ \\
\end{tabular}
\end{center}
\caption{A table of values of $\AR(K_n,C_k,q)$ for various values of $k$ and $q$, assuming $n$ is large enough.  The entries when $q=k-1$ correspond to the anti-Ramsey numbers for cycles (Theorem~\ref{thm:mbnl}).}\label{tab:cycles}
\end{table}

\section{Linear threshold for connected graphs}\label{sec:linear}

In this section we identify the linear threshold for every connected graph $G$.
We begin by showing that two standard colorings of $K_n$ yield bounds on the generalized anti-Ramsey number.

\begin{definition}[Star Coloring]\label{defn:star} Let $K_n$ be colored as follows. Totally multicolor all edges incident to a single vertex (called the TMC vertex) and then color all remaining edges with one new color. This coloring uses $n$ colors.
\end{definition}

\begin{definition}[Tur{\'a}n Coloring] Let $K_n$ be colored as follows. Totally multicolor a matching of size $\floor{n/2}$ (the TMC matching) and color all remaining edges with one new color. This coloring uses $\floor{n/2}+1$ colors.

Using these two colorings, we get the following lower bounds on the generalized anti-Ramsey number of a graph.  Recall that the matching number of a graph $G$, denoted $\nu(G)$, is the maximum size of a matching in $G$.
    
\end{definition}
\begin{proposition}\label{prop:generalDeltaNu}
    For any graph $G$ with maximum degree $\Delta(G)=\Delta$ and matching number $\nu(G)=\nu$, we have
   \begin{enumerate}[(i)]
          \item $\AR(K_n, G, \Delta+1)\ge n$ \label{p1:1}, and
       \item $\AR(K_n, G, \nu+1)\ge \floor{n/2}+1$\label{p1:2}.
   \end{enumerate} 
\end{proposition}
\begin{proof}
    To prove (\ref{p1:1}), consider a star coloring of $K_n$. Any copy of $G$ sees at most $\Delta+1$ colors. Indeed any copy of $G$ which maps vertex $v\in V(G)$ to the TMC vertex of the star coloring will see at most $d(v)+1$ colors. 

    To prove (\ref{p1:2}), consider a Tur{\'a}n coloring of $K_n$. Any copy of $G$ sees at most $\nu+1$ colors. Indeed any copy of $G$ which maps the edges of a matching of size $k$ into the TMC matching of the Tur{\'a}n coloring will see at most $k+1$ colors.
\end{proof}

 We next present two very simple lemmas.   The first guarantees the existence of some structure in graphs with given maximum degree.  Note that the \emph{broom}, denoted $B_{k,\ell}$ is the graph formed from $K_{1,\ell+1}=S_{\ell+1}$ by subdividing one of the edges $k-1$ times.  Note that Jiang and West~\cite{JW} studied the classical anti-Ramsey number of brooms.
\begin{lemma}\label{lem:hasbroom}
    Suppose $G$ is a connected graph with maximum degree $\Delta\ge 2$. Either 
    \begin{enumerate}[(i)]
        \item $|V(G)|=\Delta+1$ and 
        \begin{enumerate}
        \item $G$ is a star, or 
        \item $G$ is a triangle, or 
        \item $G$ contains a broom $B_{2,\Delta-2}$,
        \end{enumerate}
    or \item $G$ contains 
a broom $B_{2,\Delta-1}$.
    \end{enumerate}
\end{lemma}
\begin{proof}
     Consider a vertex $v$ of degree $\Delta$. If $|V(G)|=\Delta+1$, and $G$ is not a star or a triangle, then there is an edge among $N(v)$. Thus we can find a broom $B_{2,\Delta-2}$.
     If $|V(G)|>\Delta+1$, then $G$ must contain an edge from $N(v)$ to $V(G)\setminus N[v]$ since $G$ is connected. This edge plus the star centered at $v$ forms a broom $B_{2,\Delta-1}$. 
\end{proof}


The second straightforward lemma guarantees structure in graphs with given matching number.  Recall that for a graph $G$, we let $\ell G$ be the disjoint union of $\ell$ copies of $G$.

\begin{lemma}\label{lem:haspathandmatching}
    Suppose $G$ is a connected graph with matching number $\nu\ge 2$. Then $G$ contains at least one of the following:
    \begin{enumerate}
        \item $P_3 \cup (\nu-1)K_2$, or \label{lem:p3nu-1}
        \item $P_4 \cup (\nu-2)K_2$. \label{lem:p4nu-2}
    \end{enumerate}
\end{lemma}
\begin{proof}
    Let $M$ be a matching of size $\nu$ and let $e$ be an edge of $M$. Since $G$ is connected, $e$ is incident to another edge $f$. If $f\setminus e\not \in V(M)$, then we have \ref{lem:p3nu-1}. If $f\setminus e \in V(M)$, then we have \ref{lem:p4nu-2}.
\end{proof}

The following definition, introduced by Erd\H{o}s, Simonovits, and S\'os, is helpful here and elsewhere.

\begin{definition}
    Given an edge-coloring $c$ of a graph $G$, a \emph{representation} of $c$ is a subgraph of $G$ consisting of one edge of each color of $c$.
\end{definition}

Our main lemma follows.  It should be noted that we did not optimize the result in terms of the number of colors, but rather used the values that made the proof work smoothly.

\begin{lemma}\label{lem:sauer-extend}
    Let $k\ge 2$. If $E(K_n)$ is exactly colored  with $(k+2)(k+1)+1$ colors, then the coloring either contains both of
    \begin{enumerate}
        \item a rainbow $S_{k+1}$ and
        \item a rainbow $B_{2, k-1}$,
    \end{enumerate}
    or it contains all three of 
    \begin{enumerate}
        \item a rainbow $(k+1)K_2$,
         \item a rainbow $P_3\cup (k-1)K_2$, and
          \item a rainbow $P_4\cup (k-2)K_2$.
    \end{enumerate}
\end{lemma}

\begin{proof}
    Consider a representation $R$ of the coloring, that is, a graph $R$ which consists of one edge of each color. A result of Abbott, Hanson, and Sauer \cite{AHS} says that  $R$ contains either (I) a vertex of degree at least $k+2$ or (II) a matching of size $k+2$.

    If (I) happens, then we have a rainbow $S_{k+2}$ and so we certainly have a rainbow $S_{k+1}$. Suppose the rainbow $S_{k+2}$ has center vertex $v$ and leaves $a_1,\ldots, a_{k+2}$ such that edge $va_i$ receives color $i$ for all $1\le i\le k+2$. Consider the edge $a_1a_2$ in $K_n$. No matter which color $j$ this edge receives, one can find a rainbow copy of $B_{2,k-1}$ by ignoring one of $va_1$ or $va_2$ and potentially one other edge $va_j$ when $2<j\le k+2$.

    If (II) happens, then we may assume that $R$ contains a matching $M$ of size $k+2$ with edges $e_1,\ldots, e_{k+2}$ such that $e_i =v_iv_i'$ receives color $i$ for all $1\le i\le k+2$. If any edge between  $\{v_1,v_1'\}$ and $\{v_2,v_2'\}$ receives a color $j$ other than $1$ or $2$, then we have a rainbow $P_4$ which, along with $e_3,\ldots, e_{k+2}$ (potentially ignoring edge $e_j$), forms a rainbow $P_4\cup (k-1)K_2$. So without loss of generality, we may assume these 4 edges are colored with $1$ and $2$ and that $v_1v_2$ is colored with 1. So the path $v_1v_2v_2'$ and the edges $e_3,\ldots, e_{k+2}$ form a  rainbow $P_3\cup kK_2$. For the remainder of the proof, in each case, we will identify a rainbow $P_4$ and one can easily verify that there are $k-2$ edges among $e_1,e_2,\ldots, e_{k+2}$ which do not share a vertex or a color with the identified $P_4$.
    If $v_1v_3$ is colored with a color $j\not\in \{1,3\}$, then $v_1'v_1v_3v_3'$ is a rainbow $P_4$. If $v_1v_3$ is colored with $3$, then $v_2'v_2v_1v_3$ is a rainbow $P_4$. So we assume $v_1v_3$ is colored with 1. If $v_2'v_3$ is color $j\not\in \cbrac{2,3}$, then $v_2v_2'v_3v_3'$ is a rainbow $P_4$. If $v_2'v_3$ has color $3$, then $v_1v_2v_2'v_3$ is a rainbow $P_4$. So we may assume that $v_2'v_3$ is colored with $2$. Finally we consider edge $v_1v_2'$. We know this edge is colored with either $1$ or $2$. If it is colored with $1$, then $v_1v_2'v_3v_3'$ is the desired $P_4$ and if it is colored with $2$, then $v_2'v_1v_3v_3'$ is the desired $P_4$.
    \end{proof}

We are now ready to prove Theorem~\ref{thm:trees} which gives the linear threshold for connected graphs.

\begin{proof}[Proof of Theorem~\ref{thm:trees}]
    Part (\ref{thm:tree_linear}) follows directly from Proposition \ref{prop:generalDeltaNu}. Part (\ref{thm:tree_const}) follows from Lemmas \ref{lem:hasbroom}, \ref{lem:haspathandmatching}, and \ref{lem:sauer-extend}. Note that $k=\nu(G)=1$ if and only if $G$ is a triangle or a star. In this case, we need to show that every exact $(k+2)(k+1)=6$-coloring of $K_n$ contains a triangle or a star which sees at least 2 colors. This is clear (even for an exact $2$-coloring of $K_n$) since there must be some vertex where two color classes intersect. So we may assume that $G$ is not a star or a triangle and that $k\ge 2.$  
    
    First suppose that $k=\nu(G)\le \Delta(G)$. Then by Lemma \ref{lem:hasbroom}, either $|V(G)|=\Delta+1$ or $G$ contains a broom $B_{2,\Delta-1}\supseteq B_{2,k-1}$. If $|V(G)| = \Delta+1$, then since $G$ is not a star or a triangle, we must have a broom $B_{2,\Delta-2}$. This broom contains a $B_{2,k-1}$ since $k-1= \nu(G)-1\le \frac{|V(G)|}{2}-1 =\frac{\Delta+1}{2}-1 \le \Delta-2$. This last inequality holds for $\Delta\ge 3$ which we may assume since $G$ is not a triangle or a star.
    By Lemma \ref{lem:haspathandmatching} $G$ contains at least one of $P_3\cup (k-1)K_2$ or $P_4\cup(k-2)K_2$. Thus in either case of Lemma \ref{lem:sauer-extend}, we find a $(k+1)$-colored subgraph of $G$ and hence there is a copy of $G$ which sees at least $k+1$ colors.

    Similarly, if $k=\Delta(G)<\nu(G)$, then we cannot have that $|V(G)|=\Delta+1$ and so by Lemma \ref{lem:hasbroom}, $G$ contains a $B_{2,k-1}$ and also $G$ contains a $(k+1)K_2$ since we are assuming that $\nu(G)>k$. Again, in either case of Lemma \ref{lem:sauer-extend} we are done.
\end{proof}

\section{Quadratic Thresholds}\label{sec:quadratic}

Erd\H{o}s, Simonovits, and S{\'o}s (essentially) proved the following general theorem about the classical anti-Ramsey number.  Recall that given a graph $G$ and a set of graphs $\mathcal{F}$, the extremal number $\ex(G,\mathcal{F})$ is the maximum number of edges in an $\mathcal{F}$-free subgraph of $G$.

\begin{lemma}\label{lem:ESS-basic} For any graphs $G$ and $H$, with $H$ appearing as a subgraph of $G$,
    \[\ex(G, \mc{H}^-)+1\le \AR(G,H,e(H)-1) \le \ex(G,H)\]
    where $\mc{H}^- = \cbrac{H-e\,:\,e\in E(H)}$. 
\end{lemma}

This lemma is enough for us to prove part (i) of Theorem~\ref{thm:quad}.  If $G$ is bipartite on $v$ vertices with $k$ edges, then trivially $\AR(K_n,G,k)=\binom{n}2$ as there are no restrictions on the coloring.  On the other hand, by Lemma~\ref{lem:ESS-basic}, we have that
\[
    \AR(K_n,G,k-1)\leq \ex(K_n,G)=O(n^{2-2/v}),
\]
where the final estimate follows from the K\"ov\'ari-S\'os-Tur\'an theorem~\cite{KST}.

To prove part (ii) of Theorem~\ref{thm:quad}, we need to do a bit more work.
Mubayi and West proved the following lemmas which generalize Lemma~\ref{lem:ESS-basic} to $\AR(G,H,k)$ in the natural way (Lemmas 5.1 and 5.2 in \cite{MW}).  We include proofs for completeness.

\begin{lemma}\label{lem:MW-ub}  Let $\mc{F}$ be any family of subgraphs of $H$ each with at least $k+1$ edges. Further suppose that for every $F\in \mc{F}$, every copy of $F$ in $G$ can be extended to a copy of $H$ in $G$ (automatic when $G=K_n$ and $n\ge|V(H)|$). Then
    \[\AR(G, H,k) \le \ex(G, \mc{F} )\]
\end{lemma}
\begin{proof}
    Consider an exact $\ex(G, \mc{F} )+1$ coloring of $G$. Any representation of this coloring contains a $(k+1)$-colored $F\in\mc{F}$ that is contained in a copy of $H$.
\end{proof}

\begin{lemma}\label{lem:MW-lb}
    Let $\mc{F}$ be any family of graphs and suppose that the largest $\mc{F}$-free subgraph of $H$ has $k-1$ edges and that $e(G)> \ex(G,\mc{F})$. Then
    \[ \ex(G,\mc{F})+1 \le  \AR(G,H,k) \]
\end{lemma}
\begin{proof}
    Consider an $\mc{F}$-free subgraph $\hat{G}$ of $G$ with $\ex(G, \mc{F})$ many edges. Color each edge of $\hat{G}$ with a distinct color and color the edges of $G-\hat{G}$ with a new color $c$. Any copy of $H$ has at most $k-1$ edges from $\hat{G}$ and all other edges of color $c$. So every copy of $H$ sees at most $k$ colors.
\end{proof}

Notice that Lemma $\ref{lem:ESS-basic}$ follows from these two. For the upper bound, apply Lemma \ref{lem:MW-ub} with $\mc{F}=\cbrac{H}$. For the lower bound, apply Lemma \ref{lem:MW-lb} with $\mc{F}=\mc{H}^-$. In both cases, $k=e(H)-1$.

The following theorem \emph{almost} determines the quadratic threshold for all nonbipartite graphs. It is a simple corollary of Lemmas~\ref{lem:MW-ub} and~\ref{lem:MW-lb}.

\begin{theorem}
    Let $G$ be a nonbipartite graph and suppose that $k$ is the number of edges in the largest bipartite subgraph of $G$. Then $\AR(K_n, G, k-1) = o(n^2)$ and $\AR(K_n, G, k+1) = \Omega(n^2)$.
\end{theorem}
\begin{proof}
    Let $H$ be a bipartite subgraph of $G$ with $k$ edges. Then by Lemma \ref{lem:MW-ub} with $\mc{F} = \cbrac{H}$ we have
    \[\AR(K_n, G, k-1) \le \textrm{ex}(n, H)= o(n^2).\]
The last inequality follows from the Erd\H{o}s-Stone theorem~\cite{ES}.
To prove the lower bound, let $\mc{F}$ be the family of all odd cycles and apply Lemma \ref{lem:MW-lb}. The quadratic lower bound comes from K\"onig's Theorem (or really just the observation that the complete balanced bipartite graph contains no odd cycles).
\end{proof}

While one can show that $\AR(K_n,G,k)=o(n^2)$ follows from the Erd\H{o}s-Stone-Simonovits Theorem, we can do slightly better by an argument of Axenovich and K\"undgen~\cite{AK}.

\begin{theorem}
    Let $G$ be a nonbipartite graph on $v$ vertices and suppose that $k$ is the number of edges in the largest bipartite subgraph of $G$. Then  $\AR(K_n, G, k) = O(n^{2-1/(v+1)})$.
\end{theorem}
\begin{proof}
    Let $H$ be a bipartite subgraph of $G$ with $k$ edges. Let $e$ be any edge in $E(G)\setminus E(H)$. 
    Let $a = v +1$. Let $M = \ex(n,K_{a,a}) = O(n^{2-1/a})$ by the K\H{o}v\'ari-S\'os-Tur\'an theorem~\cite{KST}. Consider an exact $M+1$ coloring of $K_n$. Then any representation of this coloring contains a rainbow $K=K_{a,a}$ and say the bipartition of $K$ is $A\cup B$. Select any edge inside $A$ and suppose it has color $c$. This color $c$ appears on $K$ at most once since $K$ is rainbow. Thus we may ignore one vertex from $B$ (the one incident to color $c$ in $K$ if it exists) to obtain a rainbow $K'=K_{a,a-1}$ which does not see color $c$. Then we can embed $e$ on the edge of color $c$ and embed the rest of $H$ on $K$. This gives a copy of $H+e$, and hence $G$, which sees at least $k+1$ colors.    
\end{proof}

\section*{AI Declaration}
The mathematical content of this paper was done without the use of generative AI. Prior to submission, we used GPT 6.1 Sol to proofread.

\bibliography{genar-revised}
\bibliographystyle{amsplain}

\end{document}